\documentclass[11pt]{article}
\usepackage{amsthm,amsmath,amssymb}
\numberwithin{equation}{section}
\newtheorem{theorem}{Theorem}[section]
\newtheorem{theorem/definition}{Theorem/Definition}[section]

\newtheorem{lemma}{Lemma}[section]

\newcommand{\eqa}{\begin{eqnarray}}
\newcommand{\eeqa}{\end{eqnarray}}
\newcommand{\beq}{\begin{equation}}
\newcommand{\eeq}{\end{equation}}

\newcommand{\hb}{{\bar h}}
\newcommand{\nn}{\nonumber}
\newcommand{\p}{\partial}
\newcommand{\bu}{{\bf u}}
\newcommand{\bv}{{\bf v}}
\newcommand{\ve}{\epsilon}
\newcommand{\epf}{$\quad$\hfill
\raisebox{0.11truecm}{\fbox{}}\par\vskip0.4truecm}
\newtheorem{remark}{Remark}[section]

\newtheorem{example}{Example}[section]

\begin{document}
\title{Bihamiltonian Systems of Hydrodynamic Type and Reciprocal Transformations}
\author{Ting Xue\footnote{txue@math.mit.edu}\,,\ \ Youjin Zhang\footnote{youjin@mail.tsinghua.edu.cn}\\
{\small Department of Mathematical Sciences, Tsinghua University}\\
{\small Beijing 100084, P.R. China}}
\date{}\maketitle
\begin{abstract}
We prove that under certain linear reciprocal transformation, an evolutionary PDE of hydrodynamic type
that admits a  bihamiltonian structure is transformed to a system of the same type which
is still bihamiltonian.
\end{abstract}



\section{Introduction}
Systems of hydrodynamic type are a class of quasilinear evolutionary
PDEs of the form
\begin{equation}\label{1-2}
{\bf u}_t=V {\bf u}_x, \quad {\bf u}=(u^1,\dots, u^n)^T, \ V=(V^i_j(\bu))_{n\times n}.
\end{equation}
Assume that the above system has two
conservation laws
\begin{equation}
\frac{\partial a(u)}{\partial t}=\frac{\partial b(u)}{\partial
x},\quad\frac{\partial p(u)}{\partial t}=\frac{\partial q(u)}{\partial x}
\end{equation}
with $a(u) q(u)-p(u) b(u)\ne 0$, then we can perform a change of the
independent variables \beq\label{zh-1} (x,t)\mapsto
({y}(x,t,u(x,t)), {s}(x,t,u(x,t)) \eeq by the following defining
relations
\begin{equation}\label{eqn-51}
d{y}=a(u)dx+b(u)dt,\quad d{s}=p(u)dx+q(u)dt.
\end{equation}
Such a change of the independent variables is called a
reciprocal transformation of the system (\ref{1-2}). It originates from the
study of gas dynamics, see \cite{RS} and references therein.  Under a reciprocal transformation
the system (\ref{1-2}) remains to be a system of
hydrodynamic type
\begin{equation}\label{eqn-2}
{\bu}_s=(a\,V-b\,I)(q\, I-p\,V)^{-1} {\bu}_y\,.
\end{equation}
Here $I$ denotes the $n\times n$ unity matrix and the matrix $q I-p V$ is assumed to be nondegenerate.

The system (\ref{1-2}) is called a Hamiltonian system of hydrodynamic type if
it has the representation
\beq\label{eqn-1}
{\bu}_t=J\,\nabla h(\bu),
\eeq
where $J=(J^{ij})$ is a Hamiltonian operator of hydrodynamic type
\beq
J^{ij}=\eta^{ij}(\bu)\frac{d}{dx}-\eta^{is}(\bu)\Gamma^j_{sk}(\bu)u^k_x,\quad 1\le i, j\le n,\label{pb-1}
\eeq
the matrix $\eta=(\eta^{ij})$ is nondegenerate and symmetric on certain open subset $U$ of ${\mathbb R}^n$.
Here and henceforth summations over repeated upper and lower indices are assumed. It was proved
by Dubrovin and Novikov \cite{B.Dubrovin3, B.Dubrovin4} that  $J$ is a
Hamiltonian operator if and only if
the pseudo-Riemannian metric $(\eta_{ij})=\eta^{-1}$ is flat, and
$\Gamma^j_{sk}$ coincide with the Christoffel symbols of the
Levi-Civita connection of $(\eta_{ij})$.
So we can assume that the dependent variables $u^1,\dots, u^n$ of the system (\ref{eqn-1}) are the
{\em flat coordinates} of the metric $(\eta_{ij})$, i.e., $\eta_{ij}(u)$ are constants and $\Gamma^j_{sk}=0$.
In what follows we will also call a nondegenerate symmetric bilinear form on $T^*_{\bu} U$,
such as the the one given by $\eta=(\eta^{ij})$, a metric.

A natural question is whether under a reciprocal transformation a system of the form (\ref{eqn-1})
remains to be a Hamiltonian system of hydrodynamic type.
In the special case when the reciprocal transformation is linear in $x, t$, i.e. when
$a(u), b(u), p(u), q(u)$ are constants, Tsarev gave an affirmative answer to the above question \cite{tsarev}.
In fact,
it was shown by Pavlov \cite{Pav} that under the linear reciprocal transformation
\begin{equation}\label{eqn-4}
{y}=a x+b t, \ {s}=p\, x+q t,\quad a q-b p\ne 0 ,
\end{equation}
the Hamiltonian system (\ref{eqn-1}) with $J^{ij} = \eta^{ij} \partial_x$
is transformed to the following Hamiltonian system of hydrodynamic type:
\begin{equation}\label{eqn-3}
{\bv}_{s}={\bar J}\,\nabla \hb(\bv),\quad {\bar J}^{ij}=\eta^{ij}\,\p_y.
\end{equation}
Here the new dependent variables $v^i$ and the function $\hb(\bv)$ are defined by
\begin{eqnarray}
&&\bv=(v^1,\dots,v^n)^T=\eta \nabla (q h_0-p\, h),\quad h_0=\frac12 \eta_{ij}\, u^i u^j,\label{zh-5}\\
&&\bar{h}(\bv)=a\left(q h-\frac{1}{2} p\,
\eta^{ij}\frac{\partial h}{\partial u^i}\frac{\partial h}{\partial
u^j}\right) -b \left(q h_0-p\, (u^i\frac{\partial h}{\partial
u^i}-h)\right).
\end{eqnarray}
In the general cases, the above transformation property of a Hamiltonian system of hydrodynamic
type no longer holds true as it was shown by Ferapontov and Pavlov in \cite{FP1}.
Although under a reciprocal transformation (\ref{zh-1}), (\ref{eqn-51}) the transformed system of (\ref{eqn-1})
can still be represented as a Hamiltonian system,  the transformed
Hamiltonian operator becomes nonlocal, it contains terms with integral operator $\partial_{y}^{-1}$.

In this paper, we study the properties of a bihamiltonian system of
hydrodynamic type under the action of a reciprocal
transformation. We show that under a linear reciprocal transformation of the form (\ref{eqn-4}), a bihamiltonian
system (or a hierarchy of bihamiltonian systems) keeps to have a bihamiltonian
structure of hydrodynamic type.

The main motivation of this work comes from
the program of classification for certain class of bihamiltonian evolutionary PDEs that was
initiated by Boris Dubrovin and the second author in \cite{DZ01}.
The dispersionless limits of such evolutionary PDEs are bihamiltonian systems of hydrodynamic type.
Typical examples of
this class of evolutionary PDEs include the KdV equation and the interpolated Toda lattice equation,
their dispersionless limits will be considered in the examples of Sec.\,\ref{sec-4}.  One of the
important problems related to this classification program is whether a bihamiltonian system
of this class remains to be bihamiltonian after a linear reciprocal transformation.
The present work is a first step toward answering this problem.

The plan of the paper is as follows:
We first formulate the main results in Sec.\,\ref{sec-2}, then give their proofs
in Sec.\,\ref{sec-3}; in Sec.\,\ref{sec-4} we present two examples
to illustrate the main
results; the last section is a conclusion.

\section{The main results}\label{sec-2}
A bihamiltonian system of hydrodynamic type is a system of the form
(\ref{1-2}) that admits two compatible Hamiltonian structures of
hydrodynamic type, i.e., it has the following representation:
\beq\label{zh-2} \frac{\partial \bu}{\partial t}=J_1 \nabla
h(\bu),\quad J_1 \nabla h(\bu)\equiv J_2 \nabla f(\bu)\equiv V(\bu)
{\bu}_x. \eeq Here $J_1, J_2$ are two Hamiltonian operators of
hydrodynamic type \beq J_1^{ij}=\eta^{ij}\,\p_x,\quad
J_2^{ij}=g^{ij}(\bu)\,\p_x-g^{ik}(\bu) \Gamma^j_{kl}(\bu) u^l_x.
\eeq
The symmetric matrices $\eta=(\eta^{ij}), g=(g^{ij})$ have the properties that
$\eta$ is constant and nondegenerate, $g$ is nondegenerate on certain open
subset $U$ of $\mathbb R^n$, and $\det(g-\lambda \eta)$ does not vanish identically
for any constant parameter $\lambda$.
Compatibility of these two Hamiltonian operators means that any
linear combination $J_2-\lambda J_1$ also gives a Hamiltonian
operator of the same type. Note that we have chosen the flat
coordinates $u^1,\dots, u^n$ of the metric $\eta$ as the dependent
variables of the above system.

Let us consider the effect of the linear reciprocal transformation
(\ref{eqn-4}) on the bihamiltonian property of the system
(\ref{zh-2}). To this end, we also take into account the flow of
translation along $x$. It is also a bihamiltonian system with
respect to the above bihamiltonian structure
\beq\label{zh-4}
\frac{\p \bu}{\p t_0}=J_1 \nabla h_0(\bu),\quad J_1 \nabla
h_0(\bu)\equiv J_2 \nabla f_0(\bu)\equiv \frac{\p \bu}{\p x}.
\eeq
Here the functions $h_0, f_0$ are defined by
\beq h_0(\bu)=\frac12\,\eta_{ij}\, u^i u^j,\quad
f_0(\bu)=\frac12\, {\hat g}_{ij}\,{\hat u}^i
{\hat u}^j,\quad (\eta_{ij})=(\eta^{ij})^{-1},\ ({\hat g}_{ij})=({\hat g}^{ij})^{-1}\nn
\eeq
with ${\hat g}^{ij}$ being the
components of the metric $g$ under its flat coordinates
${\hat u}^1,\dots, {\hat u}^n$.

We introduce the new dependent variables ${\bf v}=(v^1,\dots,v^n)^T$ by the relation
(\ref{zh-5}), then the Jacobian is given by
\beq\label{zh-14}
Q:=\left(\frac{\p v^i}{\p u^j}\right)=q I-p V.
\eeq
Assume that $Q$ is nondegenerate on $U\subset {\mathbb R}^n$, denote
\beq\label{zh-14-w}
W:=Q^{-1}.
\eeq
After the reciprocal transformation (\ref{eqn-4}), the systems (\ref{zh-2}) and (\ref{zh-4}) are
transformed to
\eqa
&&\frac{\p \bv}{\p s}=(a V-b I)W {\bv}_y,\label{zh-7}\\
&&\frac{\p \bv}{\p t_0}=(a q-b p) W {\bv}_y\,. \label{zh-8}
\eeqa
Note that when the transformation (\ref{eqn-4}) satisfies the condition $a=q=0, b=p=1$, the flows $\frac{\p}{\p s}$
and $\frac{\p}{\p t_0}$ coincide.

Define two metrics ${\bar\eta}=({\bar{\eta}}^{ij}(\bv)), {\bar g}=({\bar g}^{ij}(\bv))$ whose components
in the local coordinates
$v^1,\dots, v^n$ are given by the following formulae
\begin{equation}\label{zh-9}
\bar{\eta}^{ij}(\bv)=\eta^{ij},\quad
\bar{g}^{ij}(\bv)=g^{ij}(\bu),\quad i, j =1, \dots, n,
\end{equation}
where $\bu$ and $\bv$ are related by (\ref{zh-5}).

\begin{theorem}\label{th-1}
The metrics\, ${\bar\eta},\, {\bar g}$\, are flat, the corresponding Hamiltonian operators
${\bar J}_1, {\bar J}_2$ with components
\beq
{\bar J}_1^{ij}={\bar\eta}^{ij}\, \p_y,\quad
{\bar J}_2^{ij}={\bar g}^{ij}(\bv)\,\p_y+{\bar\Gamma}^{ij}_k(\bv) v^k_y
\eeq
are compatible. Here ${\bar\Gamma}^{ij}_k(\bv)=-{\bar g}^{il} {\bar \Gamma}^j_{lk}(\bv)$ with
${\bar \Gamma}^j_{lk}(\bv)$ being the Christoffel symbols of the Levi-Civita connection of the
metric ${\bar g}$.
\end{theorem}
In other words, the metrics ${\bar\eta},\, {\bar g}$\, form a flat pencil \cite{B.Dubrovin2}.

\begin{theorem}\label{thm-1}
The systems (\ref{zh-7}) is bihamiltonian with respect to the bihamiltonian structure
${\bar J}_1, {\bar J}_2$, i.e., it has the representation
\beq
\frac{\p \bv}{\p s}={\bar J}_1 \nabla {\bar h}(\bv)\equiv {\bar J}_2 \nabla {\bar f}(\bv),
\eeq
where the functions ${\bar h},  {\bar f}$ are defined by
\beq
\frac{\p\bar h(\bv)}{\p v^i}=\frac{\p \left(a h(\bu)-b h_0(\bu)\right)}{\p u^i},
\quad \frac{\p\bar f(\bv)}{\p v^i}=\frac{\p \left(a f(\bu)-b f_0(\bu)\right)}{\p u^i},\quad
1\le i\le n.\label{zh-21}
\eeq
\end{theorem}

Now let us assume that there is given another bihamiltonian system
\beq\label{zh-10}
\frac{\p\bu}{\p t_1}=J_1 \nabla h_1(\bu),\quad J_1 \nabla h_1(\bu)\equiv J_2 \nabla f_1(\bu)\equiv A(\bu)\bu_x,
\eeq
we also assume that this flow commutes with the flow given by (\ref{zh-2}). By using Corollary 4.2 of
\cite{DLZ}, we know that the commutativity of these two flows holds true automatically
when the bihamiltonian
structure $J_1, J_2$ is semisimple, i.e. when the characteristic polynomial
$
\det(g^{ij}-\lambda\, \eta^{ij})
$
has pairwise distinct roots.

\begin{theorem}\label{th-3}
Under the reciprocal transformation (\ref{eqn-4}), the system (\ref{zh-10}) is transformed to
the form
\beq\label{zh-24}
\frac{\p\bv}{\p t_1}=(a q-b p)A W \bv_y,
\eeq
it is bihamiltonian with respect to the bihamiltonian structure ${\bar J}_1, {\bar J}_2$
and has the representation
\beq\label{zh-23}
\frac{\p \bv}{\p t_1}={\bar J}_1 \nabla {\bar h_1}(\bv)\equiv {\bar J}_2 \nabla {\bar f_1}(\bv).
\eeq
Here the functions ${\bar h}_1(\bv), {\bar f}_1(\bv)$ are defined by
\beq\label{zh-22}
\frac{\p \bar{h}_1(\bv)}{\p v^i}=(a q-b p)\frac{\p h_1(\bu)}{\p u^i},\quad
\frac{\p \bar{f}_1(\bv)}{\p v^i}=(a q-b p)\frac{\p f_1(\bu)}{\p u^i},\quad 1\le i\le n.
\eeq
\end{theorem}
\vskip 0.5cm
\begin{remark}
Associated to a bihamiltonian structure of hydrodynamic type
there is usually a hierarchy of bihamiltonian evolutionary PDEs
\beq
\frac{\p \bu}{\p t_j}=B_j(\bu) \bu_x,\quad j\ge 0.
\eeq
The independent variable $x$ can be viewed as the spatial variable, $t_j$ as the time variables, and
the flows of the hierarchy mutually commute.
From theorems \ref{th-1}, \ref{th-3} it follows that the class of bihamiltonian
hierarchies of systems of hydrodynamic type is invariant with respect to reciprocal
transformations of the form (\ref{eqn-4}).
\end{remark}

\section{Proof of the main results}\label{sec-3}

To prove the theorems of the last section, let us first prove some lemmas.
\begin{lemma} The Christoffel symbols $\Gamma^k_{ij}(\bu)$ and ${\bar\Gamma}^k_{ij}(\bv)$, expressed respectively
in the local coordinates $\bu$ and $\bv$, of the
Levi-Civita connections of the metrics $g, {\bar g}$ satisfy the following relations
\begin{equation}\label{y-4}
\Gamma^k_{ij}(\bu)=\bar{\Gamma}_{i l}^k(\bv) Q^l_j\,,\quad 1\le i, j, k\le n,
\end{equation}
where the matrix $Q$ is defined in (\ref{zh-14}).
\end{lemma}
\begin{proof}
Denote by $\ \nabla_i$ the covariant derivative along $\frac{\p}{\p u^i}$ of the Levi-Civita
connection of the metric $g$.
Then from the second Hamiltonian structure of the system (\ref{zh-2}) we get
\begin{equation}\label{c-6}
V^i_j=g^{ik}\nabla_k \nabla_j f(\bu),\quad 1\le i, j\le n.
\end{equation}
From these identities and the flatness of the metric $g$ it follows that
\beq\label{zh-12}
\nabla_k V^i_j=\nabla_j V^i_k
\eeq
for any fixed indices $1\le i,j,k\le n$. From the first Hamiltonian structure of the system (\ref{zh-2})
we also have
\beq\label{zh-17}
V^i_j=\eta^{ik}\frac{\p^2 h(\bu)}{\p u^k\p u^j},
\eeq
which lead to the identities
\beq\label{c-5}
\frac{\p V^i_j}{\p u^k}=\frac{\p V^i_k}{\p u^j}=\eta^{il} \frac{\p^3 h(\bu)}{\p u^l \p u^j \p u^k},
\eeq
so by using (\ref{zh-12}) and (\ref{zh-14}) we obtain
\beq\label{y-3}
\Gamma^i_{jl}\,V^l_k=\Gamma^i_{kl}\,V^l_j,\quad
\Gamma^i_{jl}\,Q^l_k=\Gamma^i_{kl}\,Q^l_j,\quad 1\le i, j, k\le n.
\eeq
From these relations and the identity $\nabla_k g_{ij}=0$ we obtain
\begin{eqnarray}\label{c-4}
&&\p_{u^i} g_{sj}-\p_{u^s} g_{ij}
=g_{s k}\Gamma_{ij}^k-g_{ik}\Gamma_{sj}^k
=\left(g_{s k}\Gamma_{m r}^k W^r_i -g_{i k}\Gamma_{m r}^k W^r_s\right) Q^m_j\nonumber\\
&&=\left[(\p_{u^r} g_{s m}-g_{m k} \Gamma^k_{r s}) W^r_i-
(\p_{u^r} g_{i m}-g_{m k} \Gamma^k_{r i}) W^r_s\right]Q^m_j\nonumber\\
&&=\left(\p_{u^r} g_{s m}W^r_i-\p_{u^r} g_{i m} W^r_s\right)Q^m_j
=\left(\p_{v^i} g_{s m}-\p_{v^s} g_{i m}\right)Q^m_j.
\end{eqnarray}
Substituting these identities into the formula for $\Gamma^k_{ij}$ we get
\begin{eqnarray*}
&&\Gamma^k_{ij}(\bu)=\frac{1}2 g^{ks}\left(\p_{u^j} g_{si}+\p_{u^i} g_{sj}-\p_{u^s} g_{ij}\right)\\
&&=\frac{1}2g^{ks}\left(\p_{v^m} g_{si}+\p_{v^i} g_{s m}-\p_{v^s} g_{i m}\right)Q^m_j
=\bar{\Gamma}_{il}^k(\bv) Q^l_j\,.
\end{eqnarray*}
The lemma is proved.\
\end{proof}
\begin{lemma}\label{lm-2}
The metric $\bar{g}$ is flat.
\end{lemma}
\begin{proof}
Denote by ${\bar\nabla}_i$ the covariant derivative along $\frac{\p}{\p v^i}$ of the Levi-Civita
connection of the metric ${\bar g}$, and by $R, {\bar R}$ the curvature tensors of the
metric $g, {\bar g}$ respectively
\beq
(\nabla_i\nabla_j-\nabla_j\nabla_i)\p_{u^k}={R_{ijk}}^s\, \p_{u^s},\quad
(\bar{\nabla}_i\bar{\nabla}_j-\bar{\nabla}_j\bar{\nabla}_i)\p_{v^k}
={{\bar{R}}_{ijk}}^{{}\ \ \, \,s}\, \p_{v^s}.
\eeq
By using (\ref{c-5}) and (\ref{y-4}) we get
\begin{eqnarray*}
&&{R_{ijk}}^s=\p_{u^i} \Gamma^s_{jk}-\p_{u^j} \Gamma^s_{ik}+\Gamma^s_{im} \Gamma^m_{jk}
-\Gamma^s_{jm} \Gamma^m_{ik}\\
&&=\p_{u^i}\left({\bar\Gamma}^s_{kl} Q^l_j\right)-\p_{u^j}\left({\bar\Gamma}^s_{kl} Q^l_i\right)
+{\bar\Gamma}^s_{lm} Q^l_i {\bar\Gamma}^m_{rk} Q^r_j-
{\bar\Gamma}^s_{lm} Q^l_j {\bar\Gamma}^m_{rk} Q^r_i\\
&&=\left(\p_{v^m} {\bar\Gamma}^s_{kl}-\p_{v^l} {\bar\Gamma}^s_{km}+{\bar\Gamma}^s_{mr} {\bar\Gamma}^r_{lk}
-{\bar\Gamma}^s_{lr} {\bar\Gamma}^r_{mk}\right) Q^m_i Q^l_j
={\bar{R}_{mlk}}^{\ \ \, \ \, s} Q^m_i Q^l_j,
\end{eqnarray*}
thus
$
{\bar{R}_{ijk}}^{{}\ \ \, \, s}=0
$
and we proved the lemma. \end{proof}
\noindent{\bf Proof of Theorem \ref{th-1}}\ We need to prove that
the pair of metrics $ \bar{\eta}$ and $
\bar{g}$ form a flat pencil, i.e., for any constant
$\lambda$ satisfying $\det({\bar g}^{ij}-\lambda\, {\bar\eta}^{ij})\ne 0$, the
metric ${\bar g}_\lambda={\bar g}-\lambda\, {\bar\eta}$ is flat, and the contravariant components of
its Levi-Civita connection coincide with those of ${\bar g}$ in the local coordinates $v^1,\dots, v^n$ (note
that they are flat coordinates for the metric ${\bar{\eta}}$).

In Appendix D of \cite{B.Dubrovin2}, Dubrovin gave a necessary and
sufficient condition for a pair of flat metrics to form a flat
pencil. To explain this condition, let us first recall some
notations that are introduced in \cite{B.Dubrovin2}. We denote as
above by $\nabla$ the Levi-Civita connection of the metric $g$
and by $\Gamma^{ij}_k=-g^{is} \Gamma^j_{sk}(\bu)$ its
contravariant components in the flat coordinates $u^1,\dots, u^n$
of the metric $\eta$. We also denote \beq\label{zh-15}
\p_i=\p_{u^i},\quad \p^i=\eta^{ik} \p_{u^i},\quad
\nabla_i=\nabla_{u^i},\quad \Delta^{ijk}=\eta^{is}
\Gamma^{jk}_s(\bu) .
\eeq
Then the necessary and sufficient
condition for the pair of metrics $\eta, g$ to form a flat pencil
is the existence of a vector field $\xi=\xi^i\frac{\p}{\p u^i}$,
such that the identities \beq\label{zh-18} \Delta^{ijk}=\p^i\p^j
\xi^k=\Delta^{jik},\quad g^{ij}=\p^i \xi^j+\p^j \xi^i+c^{ij} \eeq
hold true for certain constant symmetric tensor $c^{ij}$, and
\beq\label{zh-19}
\Delta_{s}^{ij}\Delta_{l}^{sk}=\Delta_{s}^{ik}\Delta_{l}^{sj},
\quad (g^{im}\eta^{jl}-\eta^{im}g^{jl})\p_m\p_l \xi^k=0,
\eeq where
$\Delta^{ij}_k=\eta_{ks}\Delta^{sij}$.

The pair of metrics $\eta, g$ related to the theorem form a flat
pencil, we are to use this fact to prove that the pair of metrics
${\bar\eta}, {\bar g}$ also fulfils the above condition. By using
the flat coordinates $v^1,\dots,v^n$ of the metric ${\bar\eta}$ we
introduce the following notations that are similar to those given in
(\ref{zh-15}): \beq\label{zh-16} {\bar\p}_i=\p_{v^i},\quad
{\bar\p}^i={\bar\eta}^{ik} \p_{v^i},\quad
{\bar\nabla}_i={\bar\nabla}_{v^i},\quad
{\bar\Delta}^{ijk}={\bar\eta}^{is} {\bar\Gamma}^{jk}_s(\bv).
\eeq
Here ${\bar\nabla}$ denotes the Levi-Civita connection of the
metric ${\bar g}$ and ${\bar\Gamma}^{ij}_k=-{\bar g}^{is}
{\bar\Gamma}^j_{sk}(\bv)$ denote the contravariant components of the
connection. From (\ref{c-6}), (\ref{zh-17}) it follows that
\beq\label{e-2-a} \eta V^T=V \eta,\quad g V^T=V g. \eeq By using the
definition (\ref{zh-14}) of $Q$  we get \beq\label{e-2} \eta Q^T=Q
\eta,\quad g Q^T=Q g.
\eeq
Then from (\ref{y-4}), the first formula
of (\ref{e-2}) and the definition (\ref{zh-9}) we get \beq
\Delta^{ijk}=\eta^{is}\Gamma^{jk}_{s}=
\eta^{is}\bar{\Gamma}^{jk}_{l} Q^l_s =\bar{\eta}^{l
s}\bar{\Gamma}^{jk}_{l} Q^i_s=\bar{\Delta}^{sjk} Q^i_s, \eeq thus by
using (\ref{zh-14}) and (\ref{zh-18}) we obtain
\begin{equation}\label{e-3}
\bar{\Delta}^{ijk}=W^i_s\Delta^{sjk}=W^i_m \p^m\p^j \xi^k
={\bar\p}^i \p^j \xi^k.
\end{equation}
Similarly, from (\ref{y-4}) and the second formula of (\ref{e-2}) we have
\beq
\Delta^{ijk}=-\eta^{is}g^{jt}\Gamma^{k}_{st}
=-\eta^{is}g^{jt}\bar{\Gamma}^{k}_{sl} Q^l_t
=-\eta^{is}g^{lt}\bar{\Gamma}^{k}_{sl} Q^j_t=\bar{\Delta}^{itk} Q^j_t,
\eeq
and consequently
\begin{equation}\label{e-4}
\bar{\Delta}^{ijk}=W^j_s\Delta^{isk}=W^j_s\p^s\p^i \xi^k
={\bar\p}^j \p^i \xi^k.
\end{equation}
Thus we arrive at the identities
\begin{equation}
{\bar\p}^i \p^j \xi^k={\bar\p}^j \p^i \xi^k,
\end{equation}
which implies the existence of a vector field $
\bar\xi=\bar\xi^i\frac{\p}{\p v^i} $ such that \beq \p^i
\xi^k={\bar\p}^i {\bar\xi}^k,\quad \le i,k \le n. \eeq So by using
(\ref{zh-18}), (\ref{e-3}) and the last formulae we obtain \beq
\bar{\Delta}^{ijk}={\bar\p}^i {\bar\p}^j {\bar\xi}^k,\quad
\bar{g}^{ij} ={\bar\p}^i {\bar\xi}^j+{\bar\p}^j {\bar\xi}^i+c^{ij}.
\eeq Now we are only left to prove the analogue of (\ref{zh-19})
for the metrics ${\bar\eta}, {\bar g}$. From (\ref{e-2-a}), (\ref{e-2}),
(\ref{e-3}) and (\ref{e-4}) we have \beq
\bar{\Delta}^{ij}_k=\bar{\eta}_{ks}\bar{\Delta}^{sij}= \eta_{ks}
W^s_t\Delta^{tij}=\Delta^{ij}_s W^s_k,\quad
\bar{\Delta}^{ij}_k=\bar{\eta}_{ks}\bar{\Delta}^{sij}
=\eta_{ks}W^i_t\Delta^{stj}=W^i_s\Delta^{sj}_k.\nn \eeq
They yield
\beq \bar{\Delta}_{s}^{ij}\bar{\Delta}_{l}^{sk}
-\bar{\Delta}_{s}^{ik}\bar{\Delta}_{l}^{sj} =W^i_m(\Delta_{s}^{m
j}\Delta_{r}^{sk}-\Delta_{s}^{m k}\Delta_{r}^{sj})W^r_l=0.\nn \eeq
Finally, by using (\ref{e-2-a}), (\ref{e-2}), (\ref{e-3}) and (\ref{e-4}) we
obtain
\begin{eqnarray*}
&&(\bar{g}^{is}\bar{\eta}^{jt}-\bar{\eta}^{is}\bar{g}^{jt})
{\bar\p}_s{\bar\p}_t {\bar\xi}^k
=\bar{g}^{is}\bar{\eta}_{l s}\bar{\Delta}^{ljk}-\bar{g}^{jt}\bar{\eta}_{l t}\bar{\Delta}^{ilk}\\
&&=g^{is}\eta_{l s} W^l_m\Delta^{mjk}-g^{jt}\eta_{l t}W^i_m\Delta^{mlk}
=g^{is}\eta_{lm} W^l_s\Delta^{mjk}-g^{jt}\eta_{l t}W^i_m\Delta^{mlk}\\
&&=g^{is} W^l_s\eta^{jm}\p_l\p_m\xi^k-
g^{jm} W^i_s \eta^{sl} \p_m\p_l\xi^k
=W^i_s(g^{sl}\eta^{jm}-\eta^{sl}g^{jm})\p_l\p_m\xi^k=0.
\end{eqnarray*}
Thus the theorem is proved. \epf

\noindent{\bf Proof of Theorem \ref{th-3}}\
Let us first prove the existence of the functions ${\bar h}_1(\bv), {\bar f}_1(\bv)$ that satisfy the equations
(\ref{zh-22}). From the commutativity condition of the flows $\frac{\p}{\p t},  \frac{\p}{\p t_1}$
\begin{equation}
\frac{\partial}{\partial t_1}(\frac{\partial \bu}{\partial
t})=\frac{\partial}{\partial t_1}(\frac{\partial \bu}{\partial t}),
\end{equation}
it follows that
\beq\label{id-2}
A V=V A,
\eeq
where the matrices $V, A$ are defined respectively in (\ref{zh-2}), (\ref{zh-10}). By using the first Hamiltonian
structure of the system (\ref{zh-10}) we obtain
\begin{equation}
\eta^{i k}\frac{\partial^2 h_{1}}{\partial u^k\partial
u^m} V^m_j=V^i_k\eta^{km}\frac{\partial^2h_{1}}{\partial
u^m\partial u^j},
\end{equation}
together with the identities given in (\ref{e-2-a}) they lead to
\begin{equation}\label{id-3}
\frac{\partial^2h_{1}}{\partial u^i\partial
u^k} W^k_j=\frac{\partial^2 h_{1}}{\partial u^j\partial u^k} W^k_i .
\end{equation}
These identities imply the existence of the function ${\bar h}_1(\bv)$
that satisfies the first equation of (\ref{zh-22}).
To prove the existence of the function ${\bar f}_1(\bv)$, let us
use the second Hamiltonian structure of the system (\ref{zh-10}) to obtain
\begin{equation}
A^i_j=g^{ik}\frac{\partial^2 f_{1}}{\partial u^k\partial u^j}
-g^{ik}\Gamma^m_{k j}\frac{\partial f_{1}}{\partial u^m}.
\end{equation}
By using (\ref{e-2-a}) and (\ref{y-3}) we have
\begin{eqnarray}
&&A^i_l V^l_j=g^{ik} V^l_j\frac{\partial^2
f_{1}}{\partial u^k\partial u^l}-g^{ik}\Gamma^m_{k l}V^l_j\frac{\partial f_1}{\partial u^m}
=g^{ik}\frac{\partial^2 f_{1}}{\partial u^k\partial
u^l} V^l_j-g^{ik}\Gamma^m_{j l} V^l_k\frac{\partial f_{1}}{\partial u^m}\nonumber\\\label{id-1}
&&=g^{ik}\frac{\partial^2 f_{1}}{\partial u^k\partial
u^l} V^l_j-g^{lk}\Gamma^m_{j l} V^i_k\frac{\partial f_{1}}{\partial u^m}.
\end{eqnarray}
Thus from the identity (\ref{id-2}) and the formula
\begin{equation}\label{id-21}
V^i_k A^k_j=V^i_k g^{k l}\frac{\partial^2 f_{1}}{\partial u^l\partial
u^j}-V^i_k g^{kl}\Gamma^m_{l j}\frac{\partial f_{1}}{\partial u^m},
\end{equation}
it follows that
\begin{equation}
V^i_k g^{k l}\frac{\partial^2 f_{1}}{\partial u^l\partial u^j}=g^{ik}\frac{\partial^2
f_{1}}{\partial u^k\partial u^l} V^l_j.
\end{equation}
So we arrive at the identities
\begin{equation}\label{id-4}
\frac{\partial^2 f_{1}}{\partial u^i\partial u^k} W^k_j=\frac{\partial^2 f_{1}}{\partial
u^j\partial u^k} W^k_i ,
\end{equation}
which imply the existence of a function ${\bar f}_1(\bv)$ that satisfies the second equation of
(\ref{zh-22}).

Now by using (\ref{zh-22}) and (\ref{y-4}) we have
\eqa &&{\bar J}_2^{ik} \frac{\p{\bar f_1}(\bv)}{\p v^k}={\bar
g}^{ik} \p_y\left(\frac{\p{\bar f}_1}{\p v^k}\right)
+{\bar\Gamma}^{ik}_l v^l_y \frac{\p{\bar f}_1}{\p v^k}\nn\\
&&=(a q-b p) \left[g^{ik} \frac{\p^2 f_1}{\p u^k \p u^m} W^m_l v^l_y
+W^m_l \Gamma^{ik}_m v^l_y \frac{\p{f}_1}{\p u^k}\right]\nn\\
&&=(a q-b p) A^i_m W^m_l v^l_y.
\eeqa
It follows that the system (\ref{zh-24}) can be represented in the form
\beq
\frac{\p\bv}{\p t_1}={\bar J}_2 \nabla {\bar f}_1(\bv).
\eeq
Similarly, one can show that the system (\ref{zh-24}) also has the expression
\beq
\frac{\p\bv}{\p t_1}={\bar J}_1 \nabla {\bar h}_1(\bv) ,
\eeq
thus it is bihamiltonian and the Theorem is proved. \epf

\noindent{\bf Proof of Theorem \ref{thm-1}}\ For the system (\ref{zh-7}) we have
\beq
\frac{\p\bv}{\p s}=\frac{\p\bv}{\p t_0} \frac{\p x}{\p s}+\frac{\p\bv}{\p t} \frac{\p t}{\p s}
=(a q-b p)^{-1} \left( a \frac{\p\bv}{\p t}-b \frac{\p\bv}{\p t_0}\right).
\eeq
The theorem follows immediately from this formula and Theorem \ref{th-3}.
\epf

\section{Two examples}\label{sec-4}
We now give two examples to illustrate the above results.
\begin{example}
For a given positive integer $m$, consider the system
\begin{equation}\label{k-1}
u_t=(m+1)\,u^m u_x.
\end{equation}
It is the $(m+1)$-th flow of the dispersionless KdV hierarchy, and has the following  bihamiltonian structure
\beq
u_t=J_1\nabla h(u)\equiv J_2 \nabla f(u) ,
\eeq
where
\beq
J_1=\partial_x,\quad J_2=u\partial_x+\frac{u_x}{2},\quad h(u)=\frac{u^{m+2}}{m+2},\
f(u)=\frac2{2 m+1}\, u^{m+1}.
\eeq
Consider the simplest linear reciprocal transformation
\begin{equation}\label{re-1}
y=t,\ s=-x.
\end{equation}
We introduce the new dependent variable $v$ as in (\ref{zh-5})  by
$
v=u^{m+1}.
$
Then after the above reciprocal transformation (\ref{k-1}) is converted to the form
\begin{equation}\label{tr-1}
v_s=-\frac{1}{m+1}\,v^{-\frac{m}{m+1}}\, v_y.
\end{equation}
By using Theorem \ref{th-1}, \ref{thm-1}, we know that
(\ref{tr-1}) is also a bihamiltonian system
\eqa
&&v_s=\bar{J}_1\nabla \bar{h}(v)\equiv \bar{J}_2 \nabla \bar{f}(v),\\
&&\bar{J}_1=\partial_y,\quad
\bar{J}_2=v^{\frac{1}{m+1}}\,\partial_y+\frac{1}{2 (m+1)} v^{-\frac{m}{m+1}}\,v_y,\\
&&
\bar{h}(v)=-\frac{m+1}{m+2}\,v^{\frac{m+2}{m+1}},\quad
\bar{f}(v)=-2 v.
\eeqa

Let us consider the $(k+1)$-th flow of the dispersionless KdV hierarchy
\beq
u_{t_1}=(k+1) u^k u_x
\eeq
which is also a bihamiltonian system
\beq
u_{t_1}=J_1\nabla h_1(u)\equiv J_2 \nabla f_1(u),\quad
h_1=\frac{u^{k+2}}{k+2},\
f_1=\frac2{2 k+1}\, u^{k+1}.
\eeq
Under the reciprocal transformation (\ref{re-1}) it is transformed to the following bihamiltonian system:
\beq
v_{t_1}=\frac{k+1}{m+1} v^{\frac{k-m}{m+1}} v_y\equiv
\bar{J}_1\nabla \bar{h}_1(v)\equiv \bar{J}_2 \nabla \bar{f}_1(v) ,
\eeq
where
\beq
{\bar h}_1(v)=\frac{m+1}{m+k+2}\, v^{\frac{k+m+2}{m+1}},\quad
{\bar f}_1(v)=\frac{2 (k+1) (m+1)}{(2 k+1)(m+k+1)}\, v^{\frac{m+k+1}{m+1}}.
\eeq
\end{example}

\begin{example}\label{exm-zh2}
Consider the long wave limit (also called dispersionless limit)
\begin{equation}\label{t-1}
u_{tt}=(e^u)_{xx}
\eeq
of the interpolated Toda equation\cite{CDZ, B.Dubrovin2}
\beq
\ve^2 u_{tt}=e^{u(x+\ve)}+e^{u(x-\ve)}-2 e^{u(x)} .
\eeq
It has the following bihamiltonian representation
\beq
\left(\begin{array}{ll}w\\ u \end{array}\right)_t=J_1\nabla
h(w,u)\equiv J_2 \nabla f(w,u) , \label{zh-21}
\eeq
where
\begin{displaymath}
J_1=\left(\begin{array}{ll}0\ &\partial_x\\ \partial_x\ &0\end{array}\right),
\ J_2=\left(\begin{array}{ll}2e^u\partial_x+e^uu_x\ &w\partial_x \\
w\partial_x+w_x\ &2\partial_x\end{array}\right),
\end{displaymath}
\begin{equation}
h(w,u)=e^u+\frac{w^2}{2},\ f(w,u)=w.\nn
\end{equation}
Let us perform the reciprocal transformation
(\ref{re-1})
and define the new dependent variables ${\bar w}, {\bar u}$ according to (\ref{zh-5})
\beq
{\bar w}=e^u,\quad \bar{u}=w.
\eeq
Then (\ref{zh-21}) is
transformed to the following equation:
\begin{equation}\label{tr-3}
\bar{w}_s=-\bar{u}_y,\ \bar{u}_s=-\frac{1}{\bar{w}}\bar{w}_y.
\end{equation}
By using Theorem \ref{th-1} and Theorem \ref{thm-1} we know that it is bihamiltonian
\begin{displaymath}
\left(\begin{array}{ll}\bar{w}\\\bar{u}\end{array}\right)_s=\bar{J}_1\nabla
\bar{h}(\bar{w},\bar{u})\equiv \bar{J}_2 \nabla
\bar{f}(\bar{w},\bar{u}),
\end{displaymath}
where
\eqa
&&\bar{J}_1=\left(\begin{array}{ll}0\ &\partial_y\\ \partial_y\
&0\end{array}\right),
\ \bar{J}_2=\left(\begin{array}{ll}2\bar{w}\partial_y+\bar{w}_y\ &\bar{u}\partial_y\\
\bar{u}\partial_y+\bar{u}_y\ &2\partial_y\end{array}\right),\\
&&\bar{h}(\bar{w},\bar{u})=-\bar{w}\log\bar{w}+\bar{w}-\frac{\bar{u}^2}{2},\\
&&\bar{f}(\bar{w},\bar{u})=-\frac{\bar{u}\log{\bar{w}}}{2}+\bar{u}-\sqrt{-4\bar{w}+\bar{u}^2}
ArcTanh(\frac{\bar{u}} {\sqrt{-4\bar{w}+\bar{u}^2}}).
\eeqa
Now we consider the flow
\begin{equation}\label{t-2}
w_{t_1}=e^uw_x+we^uu_x,\ u_{t_1}=ww_x+e^uu_x
\end{equation}
which belongs to the dispersionless Toda hierarchy\cite{CDZ}. It has the bihamiltonian
structure
\begin{displaymath}
\left(\begin{array}{ll}w\\u\end{array}\right)_{t_1}=J_1\nabla
h_1(w,u)\equiv J_2 \nabla f_1(w,u)
\end{displaymath}
with
$$
h_1(w,u)=e^uw+\frac{w^3}{6},\
f_1(w,u)=\frac{1}{2}(e^u+\frac{w^2}{2}).
$$
After reciprocal transformation (\ref{re-1}) it is transformed to
\begin{equation}\label{tr-4}
\bar{w}_{t_1}=\bar{u}\bar{w}_y+\bar{w}\bar{u}_y,\
\bar{u}_{t_1}=\bar{u}\bar{u}_y+\bar{w}_y.
\end{equation}
By using Theorem \ref{th-3}, we know that the above system
is also bihamiltonian
\begin{displaymath}
\left(\begin{array}{ll}\bar{w}\\\bar{u}\end{array}\right)_{t_1}=\bar{J}_1\nabla
\bar{h}_1(\bar{w},\bar{u})\equiv \bar{J}_2 \nabla
\bar{f}_1(\bar{w},\bar{u}),
\end{displaymath}
where
$$
\bar{h}_1(\bar{w},\bar{u})=\frac{1}{2}(\bar{u}^2\bar{w}+\bar{w}^2),\
\bar{f}_1(\bar{w},\bar{u})=\frac{\bar{u}\bar{w}}{2}.$$
\end{example}

\section{Conclusion}\label{sec-5}
We have considered the effect of a linear reciprocal transformation on the bihamiltonian structure
of a bihamiltonian system of hydrodynamic type. In the case when the bihamiltonian structure
is related to a Frobenius manifold, there are some special linear
reciprocal transformations that correspond to the Legendre transformations
between Frobenius manifolds. Such transformations
were introduced and studied by Boris Dubrovin in the setting of Frobenius manifold theory
\cite{B.Dubrovin2}.
The bihamiltonian structure and the reciprocal transformation that are considered in Example
\ref{exm-zh2} just belong to this case.

As we mentioned in the introduction, we are interested in the problem of
whether a more general class of  bihamiltonian systems
remain to be bihamiltonian after certain linear reciprocal transformations.
Such bihamiltonian systems are the objects of study of the classification program that was initiated in
\cite{DZ01}, their dispersionless limits are bihamiltonian systems of hydrodynamic type.
A positive answer to this problem was given
in \cite{CDZ} for the example of the extended Toda lattice hierarchy, a hierarchy of bihamiltonian evolutionary PDEs.
It was shown there that after certain
linear reciprocal transformation this hierarchy is transformed to the extended nonlinear Schr\"odinger
hierarchy which is also bihamiltonian. Results regarding the general case for this problem
will appear in \cite{DZ05} and other publications.

\vskip 0.2truecm \noindent{\bf Acknowledgments.} The authors are
grateful to Boris Dubrovin for encouragements and helpful discussions.
The researches of Y.Z. were partially
supported by the Chinese National Science Fund for Distinguished
Young Scholars grant No.10025101 and the Special Funds of Chinese
Major Basic Research Project ``Nonlinear Sciences''.

\end{document}